\documentclass[twoside,11pt]{article} 
\usepackage{amssymb}
\usepackage{graphicx}^^M
\usepackage{amsthm}
\newfont{\bbb}{msbm10 scaled\magstep 1}
\newcommand{\inter}{{\rm int}}
\newcommand{\cen}{{\rm cen_y}}
\newcommand{\area}{{\rm area}}
\font\bbigbold=cmbx10 at 17 pt 
\date{}

\title
{\bbigbold 
\noindent 
Position of the centroid of a planar convex body}

\begin{document}

\baselineskip 18pt 

\maketitle

\noindent
\author{MAREK LASSAK}

\vskip 0.2cm
\pagestyle{myheadings} \markboth{Marek Lassak}{Centroid of a planar convex body}

\vskip 0.2cm
\noindent
{\bf Abstract.}
It is well known that any planar convex body $A$ permits to inscribe an affine-regular hexagon $H_A$. We prove that the centroid of $A$ belongs to the homothetic image of $H_A$ with ratio $\frac{4}{21}$ and the center in the center of $H_A$. This ratio cannot be decreased.



\date{}

\maketitle

\section{Introduction}
This paper concerns the position of the centroid of a planar convex body, i.e., a closed bounded convex set.
Recall that the notion of centroid is discussed by, among others, Bonnesen and Fenchel \cite{[B+F]}, Gr\" unbaum \cite {[G]}, Hammer \cite{[H]} and Neumann \cite{[N]}.

As usual, by an {\it affine-regular hexagon} we understand a non-degenerated affine image of the regular hexagon.
Besicovitch \cite{[B]} proved that for every planar convex body $A$ there exists an affine-regular hexagon $H_A$ inscribed in $A$.
Our aim is to prove that the centroid of $A$ belongs to the homothetic image $\frac{4}{21}H_A$ of $H_A$ with ratio $\frac{4}{21}$ and the center in the center of $H_A$.
In general, this ratio cannot be lessened, which is explained at the end of the paper.

For a compact set $C$ of the Euclidean plane $E^2$ denote by ${\rm cen}_x (C)$ and $\cen (C)$ the first and the second coordinates of the centroid of $C$.
Let compact sets $B_1, \dots , B_n \subset E^2$ with non-empty interiors have disjoint interiors and $B = \bigcup_{j=1}^n B_j$.
It is well known that 

$${\rm cen}_x(B) = \frac{\Sigma_{j=1}^n {\rm cen}_x(B_j) \cdot{\rm area}(B_j)}{\Sigma_{j=1}^n  \area (B_j)}, 
\ \ \ 
{\rm cen}_y(B) = \frac{\Sigma_{j=1}^n {\rm cen}_y(B_j) \cdot{\rm area}(B_j)}{\Sigma_{j=1}^n  {\rm area}(B_j)}. \eqno (1)$$

\section{The position of the centroid of a convex body with respect to an inscribed affine-regular hexagon}

Let $D \subset E^2$ and $\ell$ be a straight line.
Imagine $D$ as the union of segments (including one-point segments) being intersections of $D$ by straight lines perpendicular to $\ell$.
Shift every such a segment perpendicularly to $\ell$ in order to obtain its image centered at $\ell$. 
Denote the union of all these obtained segments by ${\rm sym}_\ell D$. 
It is the result of the Steiner symmetrization of $D$.
The proof of the following lemma is given in a number of books. 
For instance in Section 40 of \cite{[B+F]}.

\vskip0.25cm
\noindent
{\bf Lemma.} 
{\it If $D \subset E^2$ is convex, then ${\rm sym}_\ell D$ is convex.}

\vskip0.2cm
\noindent
{\bf Theorem.}
{\it Let $A \subset E^2$ be a convex body and $H_A$ be an affine-regular hexagon inscribed in $A$. 
Then the centroid of $A$ belongs to the homothetic image of $H_A$ with ratio $\frac{4}{21}$ and center in the center of $H_A$.}

\begin{proof}
For better clarity, we divide the proof into a preliminary text mostly on notations, and then Parts 1--8 with considerations.

We do not lose the generality assuming that the successive vertices $a_1, \dots , a_6$ of $H_A$ are $(1,1)$, $(-1,1)$, $(-2,0)$, $(-1,-1)$, $(1,-1)$, $(2,0)$, see Figure. 
Denote by $o$ the center $(0,0)$ and by $a$ the midpoint of $a_1a_2$.
Since we deal with ${\cen}(A)$, by Lemma we may assume that $x=0$ is an axis of symmetry of $A$.

In order to prove the assertion, let us show that for any side of $\frac{4}{21}H_A$ the centroid of $A$ is on the same side of the straight line containing this side which contains $o$.
Observe that it is enough to show this for one side of the hexagon $\frac{4}{21}H_A$. 
Let us provide this task for the side connecting $\frac{4}{21}a_1$ and $\frac{4}{21}a_2$.   

Denote by $\overline{a}_i$ the intersection of the straight lines containing $a_ia_{i+1}$ and $a_{i-1}a_{i-2}$ for $i=1,\dots,6$ (mod $6$), see Figure. 
We define the star $S(H_{A})$ over $H_A$ as the union of $H_A$ and six triangles $T_i(H_A) = a_{i-1}{\buildrel - \over a_i}a_i$, where $i=1, \dots , 6$ and where $a_0$ means $a_6$.
From the convexity of $A$ we conclude that $A \subset S(H_A)$.

We do not make our considerations narrower assuming that the centroid of $A$ is over or on the axis $y=0$.
Since our aim is to show that ${\cen}(A) \leq \frac{4}{21}$ for every convex body $A$, it is sufficient to consider only such convex bodies $A$ which are disjoint with the interiors of $T_4(H_A)$, $T_5(H_A)$ and $T_6(H_A)$.
Still the closure of $A \setminus \bigcup_{i=4}^6 T_i(H_A)$ is a convex body with $H$ inscribed and the centroid at the same or higher level.

Provide any supporting straight line $L_1$ of $A$ at $a_1$ and the symmetric (with respect to $x=0$) supporting line $L_2$ of $A$ at $a_2$.
Denote by $u=(0,w)$ the intersection point of $L_1$ (and thus of $L_2$) with the axis $x=0$. 
Since the second coordinates of $a$ and $\overline{a}_2$ are equal to $1$ and $2$, respectively, we have $w \in [1,2]$.

Since $L_1$ passes through $u=(0,w)$ and $a_1 =(1,1)$, it has the equation 
$y-1 = (-w+1)(x-1)$.
Its point of intersection with the segment $a_6{\buildrel - \over a_1}$ (being a subset of the straight line $y=x-2$) is $m_1 =(\frac{2+w}{w}, \frac{2-w}{w})$.
Similarly, we get the symmetric point $m_2$ being the intersection of $L_2$ with the segment $a_3{\buildrel - \over a_3}$.

\newpage

\begin{figure}[htbp]
\centering
\hskip0.2cm
\includegraphics[width=15.4cm,height=8.1cm]{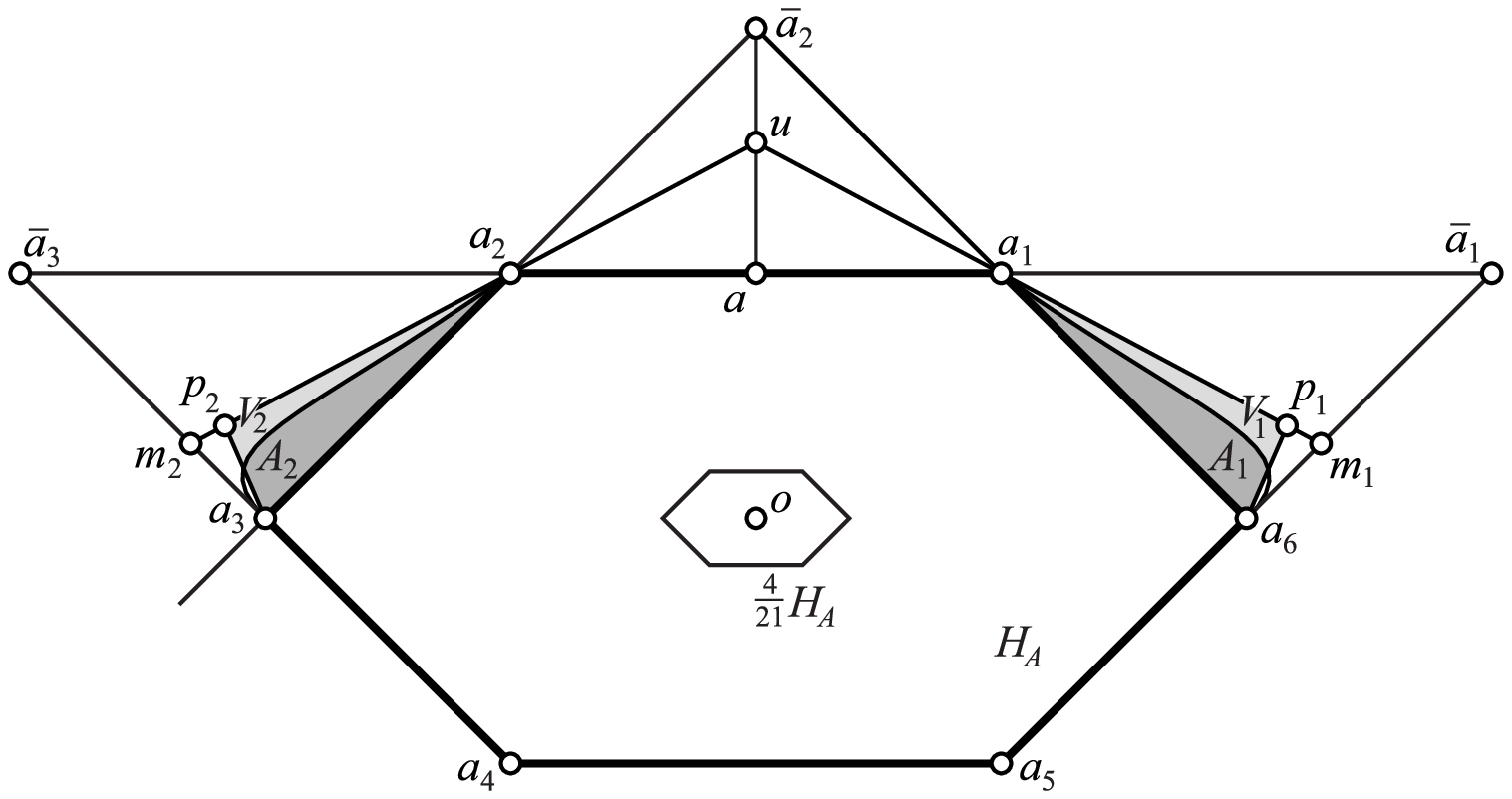} \\ 
\vskip0.2cm
\centerline
{Figure. Illustration to the proof of Theorem}
\end{figure}

Later we explain the geometric meaning of the following number 

$$w_0 = \frac{1}{3}({\root 3 \of {44-3\sqrt {177}}} + {\root 3 \of {44 +3\sqrt {177}}} -1) = 1.6589670... .$$ 

Parts 3--7 lead to the proof of our theorem for $w \in [w_0, 2]$ and Part 8 for $w \in [1, w_0]$ .

\bigskip
{\bf Part 1} where we introduce a heptagon and find its $\cen$.

Let $z \in [0,1]$.
Since $a_1 =(1,1)$ and $m_1 =(\frac{2+w}{w}, \frac{2-w}{w})$, every point $p_1(z)$, or shortly $p_1$, of $a_1m_1$ has the form $(1-z)a_1 + zm_1$.
So $p_1 = ((1-z) +z\frac{2+w}{w}, (1-z) +z\frac{2-w}{w}) = (z\cdot \frac{2}{w} +1, z\cdot\frac{2-2w}{w}+1)$.   
The symmetric point with respect to $x=0$ is denoted by $p_2$.
The second coordinates of them are $z \cdot \frac{2-2w}{w} +1$.

Consider the heptagon $G= up_2a_3a_4a_5a_6p_1$.
The area of each of the two symmetric {\it wings} $W_1 = a_1a_6p_1$ and $W_2 = a_2a_3p_2$ of $G$ is $z\cdot \frac{2-w}{w}$ and $\cen$ of each wing of this heptagon is $\frac{1+ z\frac{2-2w}{w}+1}{3} = \frac{2+ z\frac{2-2w}{w}}{3}$.
The area of the triangle $a_1ua_2$ is $\frac{1}{2}\cdot 2\cdot (w-1) = w-1$ and its $\cen$ is $\frac{2+w}{3}$.
Moreover, the area of $H_A$ is $6$ and its $\cen$ is $0$.
Taking all this into account and having in mind that $G = H_A \cup a_1a_6p_1 \cup a_2a_3p_2 \cup a_1ua_2$,
 by the right part of (1) we conclude that  

$$\cen (G) = \frac{0+\frac{2}{3}(2+z\frac{2-2w}{w})z\frac{2-w}{w} + \frac{2+w}{3}(w-1)}{6+ 2z\frac{2-w}{w} +w-1}$$

\noindent
which, after a simplification, equals to
 
$$\frac{2(2+z\frac{2-2w}{w})z\frac{2-w}{w} + w^2+w-2} {6z\frac{2-w}{w}+3w+15}. \hskip 0.5cm \eqno (2)$$

\bigskip
{\bf Part 2}  whose aim is to show the following statement

\medskip
{\it Denote by $\nu$ the numerator and by $\delta$ the denominator of $\cen (G)$ as in {\rm (2)} (so $\cen (G) = \frac{\nu}{\delta}$).
Consider a truncation of the wings $W_i$ of $G$ to symmetric convex subsets $A_i = W_i \cap A$ for $i=1,2$.
Put $V_i= W_i \setminus A_i$ for $i=1,2$ and $V= V_1 \cup V_2$.
We have}

$$\frac {\nu - \area (V) \cen(V)}{\delta- \area (V)} \leq 
\frac{\nu}{\delta} 
\hskip 0.4cm {\it iff}  \hskip 0.4cm \cen (V)  \geq 
\frac{\nu}{\delta}. \eqno (3)
$$

\medskip
Let us confirm this.
We have
$\frac {\nu - \area (V) \cen(V)}{\delta- \area (V)} \leq \frac{\nu}{\delta}$
iff
$\delta (\nu - \area (V) \cen(V)) \leq \nu(\delta- \area (V))$ 
iff 
$\nu \cdot\area (V) \leq \delta \cdot \area (V) \cen(V)$ 
iff 
$\nu \leq \delta \cdot \cen(V)$ 
iff $\cen (V)  \geq \frac{\nu}{\delta}$.

\medskip
Observe that $\frac {\nu - \area (V) \cen(V)}{\delta- \area (V)}$ is nothing else but $\cen (A')$, where $A' = G \setminus V$.

\bigskip
{\bf Part 3} where we start considerations for $w \in [w_0, 2]$.  
 
\medskip
For every $w \in [w_0,2]$ we are looking for the positions of $p_1$ and thus of $p_2$ such that $\cen (G)$ is the largest.
For this reason let us find the derivative of the function $(2)$ with respect to $z$:

$$\frac{2(w-2)[4z^2(-w^2+3w-2) +4z(w^2+4w-5) + (w^4-w^3-12w^2)]}{3w(w^2-2wz+5w+4z)^2}. \eqno (4)$$

The discriminant of the quadratic function in the square bracket is $16w^2(2w^4+ 4w^3-w^2-6w+1)$.
Hence (4) equals $0$ for $z = \frac{w(w^2+4w-5 \pm \sqrt{2w^4+ 4w^3-w^2-6w+1})}{2(w^2-3w+2)}$.
Take into account only the root

$$z_w = \frac{w(w^2+4w - 5 - \sqrt{2w^4+ 4w^3-w^2-6w+1})}{2(w^2-3w+2)} \eqno (5)$$

\noindent
which is positive for every $w \in [w_0,2)$ (the other one is always negative here).
Moreover, put $z_2 = \lim_{w \to 2^-}z_w$.
This is $z_2 = \frac{5}{7}$. 

We see that for any fixed $w \in [w_0, 2]$ the global maximum of $(2)$ as a function of $z$ from the interval $[0,1]$ can be only for $z=0$, $z= z_w$ or $z=1$.
Substituting these three $z$ into $(2)$ we see that the global maximum of (2) in the interval $[0,1]$ is at $z=z_w$ for every fixed $w \in [w_0, 2]$.

\bigskip
{\bf Part 4} where our aim is to show that for each $w \in [w_0,2]$ the value of (2) for $z=z_w$ is at most $\frac{4}{21}$.

\medskip
This task with substituting $z=z_w$ into (2) seems to be very complicated to perform.
We can get it around by performing the more general task to show that for every $w \in [w_0 ,2]$ and $z \in [\frac{5}{7},1]$ we have

$$\frac{2(2+z\frac{2-2w}{w}) z\frac{2-w}{w} +w^2+w-2}{6z\frac{2-w}{w}+3w +15} \leq \frac{4}{21}. \eqno(6)$$

This task is more general since $z_w$ belongs to $[\frac{5}{7},1]$ for every $w \in [w_0, 2]$.
Really, the inequality $z_w \leq 1$ is equivalent to $w^6 -7w^4 -2w^2 +16w^2 +8w -16 \geq 0$ and thus to $(w-1)(w-2)(w+3)(w^3 +w^2 -2w -4) \geq 0$, which means that it holds true in $[1,2]$ if and only if $w \in [w_0, 2]$ (still $w_0$ is the only real root of this polynomial).
Moreover, the inequality $\frac{5}{7} \leq z_w$ is equivalent to $-7(w-2)(7w^2+22w+20) \geq 0$, which means that it holds true in the whole interval $[1,2]$, so in particular for every $w \in [w_0, 2]$.

Equivalently to (6), it is sufficient to show that 

$$28z^2(w^2-3w+2) +20zw(2-w) +w^2(7w^2+3w-34)  \eqno(7)$$

\noindent
is at most $0$ for every point $(w,z)$ of the rectangle $[w_0,2] \times [\frac{5}{7}, 1]$.

In order to simplify evaluations consider this task in the larger rectangle $[1,2] \times [\frac{5}{7}, 1]$.

Let us apply the following  method of finding the global maximum of a continuous function $f(w,z)$ in a polygon $R \subset E^2$.
Namely, first we find the points being the solutions of the system of two equations when partial derivatives of our function $f(w,z)$ are $0$ in the interior of $R$.  
Next we write the equations of the sides in the forms $z=g(w)$ or $w=g(z)$.
We find the critical points in the relative interiors of each side, where the derivative of the respective equation is $0$.
Finally, we check the values of $f(w,z)$ at the vertices of $R$. 
The largest value at all the found points gives the maximum value of $f(w,z)$ in $R$. 

In our particular case our function $f(w,z)= 28z^2(w^2-3w+2) +20zw(2-w) +w^2(7w^2+3w-34)$ is given by (7). 
Moreover, $R= [1,2] \times [\frac{5}{7}, 1]$.
According to the recalled method we find the partial derivatives 
$f'_w (w,z) = 28w^3 +9w^2 +56wz^2 -40wz -68w -84z^2 +40z$ and $f'_z (w,z)  = 56w^2z -20w^2 -168wz +40w +112z$.
Consider the system of equations when both are $0$.
Finding $z$ from the second and substituting to the first we get three solutions: $w \approx -1.8$, $w=0$ and $w \approx 1.544$.
None of them is in the interval $[w_0,2]$.
Hence the system of equations has no solution in our $R$, and thus in its interior.

Let us find the critical points in the relative interiors of the sides.
After substituting $z= \frac{5}{7}$ to $f(w,z)$ we get $\frac{1}{7} (49w^4 +21w^3 -238w^2 -100w +200)$.
Its derivative $\frac{1}{7} (196w^3 +63w^2 -476w -100)$ is $0$ only at $w_1 =1.5103...$.
Placing $z= 1$ to our $f(w,z)$ we get $7w^4 +3w^3 -26w^2 -44w +56$.
Its derivative $28w^3+9w^2- 52w -44$ equals $0$ in $[1,2]$ only at $w_2 = 1.5427...$.
Substituting $w= 1$ to $f(w,z)$ we get $20z-24$, which is negative for every $z \in [\frac{5}{7},1]$.
Placing $w= 2$ to $f(w,z)$ we get $0$ for every $z \in [\frac{5}{7},1]$.

We have 
$f(w_1, \frac{5}{7}) \approx -23.803$,
$f(w_2, 1) \approx -23.094$,
$f(1, 1) =-4$, 
$f(1, \frac{5}{7}) =-9.714$, 
$f(2, \frac{5}{7})= 0$, 
and $f(2,1) =0$.
Thus the global maximum of $f(w,z)$ in $R$ is $0$.
Hence (7) is at most $0$ and thus (6) holds true in $R$. 
We conclude that (2) for $z=z_w$ is at most $\frac{4}{21}$ for every $w \in [1,2]$ and so for every $w \in [w_0, 2]$.

\bigskip
{\bf Part 5} where we show that $\cen (V) \geq \frac{\nu}{\delta}$ for any $w \in [w_0, 2]$ and $z_{w}$ in place of $z$.  

\medskip
Recall from Part 2 that $\frac{\nu}{\delta} = \cen (G)$.
Looking at the second coordinates of $a_1, a_6$ and $p_1$ we get $\cen (a_1a_6p_1) = (z_w\frac{2-2w}{w} +1)/2$. 
Hence
$\cen(V_1) \geq (z_w\frac{2-2w}{w} +1)/2$ and so $\cen(V) \geq (z_w\frac{2-2w}{w} +1)/2$.

We see that in order to confirm the promise of Part 5 it is sufficient to show that   

$$\frac{z_w\frac{2-2w}{w} +1}{2} \leq  \frac{2(2+z_w\frac{2-2w}{w})z_w\frac{2-w}{w} + w^2+w-2}{6z_w\frac{2-w}{w}+3w+15}, \eqno(8)$$ 

\noindent
for every $w \in [w_0,2]$, where the right side is taken from (2).
Instead, let us show the inequality 

$$\frac{z\frac{2-2w}{w} +1}{2} \leq  \frac{2(2+z\frac{2-2w}{w})z\frac{2-w}{w} + w^2+w-2}{6z\frac{2-w}{w}+3w+15},$$ 

\noindent
or equivalently, let us show that

$$8z^2 -12z^2w -22zw^2 +26zw +4z^2w^2 -6zw^3 -2w^4 +w^3 +19w^2 \eqno(9)$$

\noindent
is at most $0$ for every point $(w,z)$ of the piece of the curve $z=z_w$ when $w \in [w_0,2]$.

Instead, let us find the global maximum of (9) in a triangle containing it.
Namely, in the triangle $T$ between the straight lines $w=2$, $z= -\frac{5}{7}w +\frac{15}{7}$, and $z=1$. 
Its vertices are $(2,1)$, $(\frac{8}{5}, 1)$ and $(2, \frac{13}{21})$.

First let us show that the piece of the curve
$z=z_w$ for $w \in [w_0,2]$ is a subset of $T$. 
The reason is that 
$-\frac{5}{7}w +\frac{15}{7} \leq z_w \leq 1$ for every $w \in [w_0,2]$.
The left inequality is equivalent to the inequality 
$(2839w^4 -10571w^3 +18960w^2 -284088w +22472)(w-1)(2-w) \geq 0$ which holds true for every $w \in (-\infty, \infty)$.
Thus in $[1,2]$ and so for every $w \in [\frac{8}{5},2]$. 
The right inequality $z_w \leq 1$ is shown just after (6).

Next let us find the global maximum of (9) in $T$ by the method described in Part 4.

Consider the system of equations 
$-8w^3 -18w^2z +3w^2 +8wz^2 -44wz +38w -12z^2+26z=0$ and 
$-6w^3 +8w^2z -22w^2 -24wz +26w +16z=0$ (where the left sides are the partial derivatives of (9)).
Finding  
$z = \frac{3w^3 +11w^2 -13w}{4w^2 -12w +8}$ 
from the second and substituting it into the first we get the equation 
$w(68w^6 -141w^5 -262w^4 +359w^3 +356w^2 -447w +68)=0$ whose solutions are $w=0, w \approx 0.183, w \approx 0.951, w \approx 1.037$ and $w \approx 2,614$.
All these $w$ are out of the interval $[\frac{8}{5}, 2]$ which implies that
all the obtained points $(w,z)$ are out of $T$.
Thus the system of equations has no solution in the interior of $T$. 

Look for critical points in the relative interiors of the sides.
Substituting $z= -\frac{5}{7}w +\frac{15}{7}$ into (9) we get $\frac{1}{49} (212w^3 -511w^2 -789w +1830)$. 
Its derivative $\frac{1}{49} (848w^3 -1533w^2 -1578w +1830)$ is never $0$ in $[\frac{8}{5},2]$.
Putting $z=1$ into (9) we get $-2w^4 -5w^3 +w^2 +22w$.
Its derivative $-8w^3 -15w^2 +2w +22$ is never $0$ in $[\frac{8}{5}, 2]$. 
Putting $w=2$ into (9) we get $8z^2-84z+52$. 
Its derivative $16z-84$ is never $0$ in $[\frac{5}{7},1]$.

The value of (9) at $(2, 1)$ is $-44$, at $(\frac{8}{5}, 1)$ is $-11.827...$, and at $(2, \frac{13}{21})$ is $-4.598...$.
So the global maximum of the function (9) in $T$ is $-11.827...$.
Hence (9) is always negative in $T$.

Consequently, we have shown that (9) is at most $0$ in $T$ and thus that (8) is true for every $w \in [w_0,2]$.
Therefore $\cen (V) \geq  \frac{\nu}{\delta}$ for $G$ with $z_w$ in the part of $z$.

\bigskip
{\bf Part 6} where we show that $\cen (A') \leq \frac{4}{21}$ for $w \in [w_0,2]$.

\medskip
Recall that $\cen (G) = \frac{\nu}{\delta}$.
By Part 5 and by (3) we have 
$\frac {\nu - \area (V) \cen(V)}{\delta- \area (V)} \leq \frac{\nu}{\delta}$.
The left side is $\cen (A')$ and the right one is $\cen (G)$ with $G$ is taken for $z=z_w$.
By (6) it is at most $\frac{4}{21}$.
So $\cen (A') \leq \frac{4}{21}$.

\bigskip
{\bf Part 7} on enlarging $A'$ up to $A$ which leads to the proof of our theorem for $w \in [w_0,2]$.

\medskip
Put $A_1'' = A \cap p_1a_6m_1$, $A_2'' = A \cap p_2a_3m_2$ and $A'' = A_1'' \cup A_2''$.
Clearly $A = A' \cup A''$.
We intend to show that adding $A''$ to $A'$ does not increase $\cen$, so that $\cen (A) \leq \cen (A')$. 

First let us show that if the triangles $p_1a_6m_1$ and $p_2a_3m_2$ are added to $A'$, then $\cen$ does not increase.
Applying the easy to show implication: ``if $\inter(X) \cap \inter (Y) = \emptyset$, $\cen (X) \leq \mu$ and  $\cen (Y) \leq \mu$, then $\cen (X \cup Y) \leq \mu$ as well" and having in mind that $\cen(A') \leq \frac{4}{21}$ (see Part 6), it is sufficient to show that $\cen (p_1a_6m_1) \leq \frac{4}{21}$ (then also $\cen (p_2a_3m_2) \leq \frac{4}{21}$). 
Let us show this.
Since  $\cen (p_1a_6m_1) = (z_w \cdot \frac{2-2w}{w} +1 + \frac{2-w}{w})/3$, we have to show that this is at most $\frac{4}{21}$. 
This task is equivalent to $7z_w(2-2w) \leq 4w -14$.
After substituting $z_w$ and providing some simplifications, this inequality is equivalent to $h(w) \geq 0$, where $h(w)= 49w^6 -196w^5 +105w^4 +1946w^2 +1800w -400$. 
We have $h''(w) = 14(105w^4 -280w^3 +90w^2 +278)$.
From the fact that $h''(w)$ always positive we conclude that $h'(w)$ is an increasing function.
Thus from $h(1) = 3304$ we see that $h(w) \geq 0$ for $w\geq 1$, and thus for every $w \in [w_0, 2]$. 
Hence $\cen (p_1a_6m_1) \leq \frac{4}{21}$.

Also for adding only $A''$ to $A'$, the value of $\cen$ does not increase.
The reason is that $\cen (A_1'') \leq \cen (p_1a_6m_1)$ and 
$\cen (A_2'') \leq \cen (p_2a_3m_2)$.
The first follows from the convexity of $A_1''$ and from the observation that every segment jointing $a_6$ with a point of $p_1m_1$ has in common with $A_1''$ only a segment which is lower.
Analogously, we confirm the second inequality.

We conclude that  $\cen (A) \leq \cen (A')$.
This and $\cen (A') \leq \frac{4}{21}$ (see Part 6) imply $\cen (A) \leq \frac{4}{21}$.

\bigskip
{\bf Part 8} where we prove our theorem for $w \in [1, w_0]$.

\medskip
Consider the pentagon $P= m_1um_2a_4a_5$.
It is the special case of $G$ for $z=1$.
Thus substituting $z=1$ to (2) we see that $\cen (P)$ equals to $\frac{w^4+w^3-2w^2-4w+8}{3w(w^2+3w+4)}$.
In order to show that this is at most $\frac{4}{21}$ for every $w \in [1, w_0]$
take into account the equivalent inequality $(w-2)(7w^3 +17w^2 +8w -28) \leq 0$.
Its left side equals $0$ only for $w=2$ and $w = w_3 \approx 0.934$. 
Consequently, this inequality and thus also the preceding one hold true in $[w_3,2]$.
Hence also in $[1,w_0]$.
Resuming, $\cen (P) \leq \frac{4}{21}$.

From $z=1$ we see that $V_1 = a_1a_6m_1$ for our $P=G$.
The second coordinates of $a_1, a_6$ and $m_1$ give $\cen (a_1a_6m_1)=  \frac{4-2w}{3w}$. 
Hence $\cen (V_1) \leq \frac{4-2w}{3w}$.
Thus by $\cen (V) = \cen (V_1)$ we get $\cen (V) \leq \frac{4-2w}{3w}$.

In order to show that the right side of (3) is now true we have to show that $\cen (V) \geq \cen (P)$, where, as $P$ takes the role of $G$, the right side is denoted by $\frac{\nu}{\delta}$ in (3).
Hence we have to show that 
$\frac{4-2w}{3w} \geq \frac{w^4+w^3-2w^2-4w+8}{3w(w^2+3w+4)}$.
This is equivalent to the inequality
$w^4 +3w^3 -8w -8 \leq 0$.  
A simple evaluation confirms that it is true in $[1,w_0]$ (by the way, we have here the equality just for $w=w_0$).
Thus $\cen (V) \geq \cen (P)$.

The shown inequality means that the right side of (3) is fulfilled.
Hence the left side of (3), i.e., $\cen (A') \leq \frac{\nu}{\delta}$ holds true. 
From $A = A'$ for our $P=G$ we obtain $\cen(A) \leq \frac{\nu}{\delta}$.
Consequently, from $\frac{\nu}{\delta}  = \cen (P) \leq \frac{4}{21}$ we conclude that $\cen(A) \leq \frac{4}{21}$.

\medskip
Thanks to results of Parts 7 and 8 the thesis of our theorem holds true.
\end{proof}

The ratio $\frac{4}{21}$ in Theorem cannot be lessened as it follows from the example of the pentagon ${\buildrel - \over a_2}a_3a_4a_5a_6$ in the part of $A$, and the hexagon $a_1 \dots a_6$ as $H_A$.
The author expects that there are no more such examples besides the affine images of the above presented one.

\vskip0.15cm
\noindent
Marek Lassak

\noindent
University of Science and Technology

\noindent
85-789 Bydgoszcz 

\noindent
Poland

\noindent
e-mail: lassak@pbs.edu.pl

\end{document}